\RequirePackage{ifpdf}
\ifpdf 
\documentclass[pdftex]{sigma}
\else
\documentclass{sigma}
\fi

\newcommand{\R}{\mathbb R}
\newcommand{\E}{\mathcal E}

\newcommand{\g}{\mathfrak g}

\newcommand{\gl}{\mathfrak{gl}}
\newcommand{\gh}{\mathfrak h}

\newcommand{\m}{\mathfrak m}
\newcommand{\gsl}{\mathfrak sl}

\newcommand{\gr}{\operatorname{gr}}
\newcommand{\Hom}{\operatorname{Hom}}
\newcommand{\im}{\operatorname{Im}}
\renewcommand{\ker}{\operatorname{Ker}}
\newcommand{\tr}{\operatorname{tr}}

\newcommand{\dd}{\operatorname{d}}
\newcommand{\dx}{\operatorname{dx}}

\begin{document}

\allowdisplaybreaks

\renewcommand{\PaperNumber}{076}

\FirstPageHeading

\ArticleName{Third Order ODEs Systems\\ and Its Characteristic Connections}

\ShortArticleName{Third Order ODEs Systems and Its Characteristic Connections}

\Author{Alexandr MEDVEDEV}

\AuthorNameForHeading{A.~Medvedev}

\Address{Faculty of Applied Mathematics, Belarusian State University, \\  4, Nezavisimosti
Ave., 220030, Minsk, Republic of Belarus}

\Email{\href{mailto:Sasha.Medvedev@gmail.com}{Sasha.Medvedev@gmail.com}}

\ArticleDates{Received April 20, 2011, in f\/inal form July 27, 2011;  Published online August 03, 2011}

\Abstract{We compute the characteristic Cartan connection associated with a system of third order ODEs. Our connection is dif\/ferent from Tanaka normal one, but still is uniquely associated with the system of third order ODEs. This allows us to f\/ind all fundamental invariants of a system of third order ODEs and, in particular, determine when a system of third order ODEs is trivializable. As application dif\/ferential invariants of equations on circles in $\R^n$ are computed.}


\Classification{34A26; 53B15}

\section{Introduction}
\subsection{Dif\/ferential equation as a structure on a f\/iltered manifold}

The main purpose of this article is to study geometry of systems of ordinary dif\/ferential equations of third order. The geometry of ordinary dif\/ferential equations or, more generally, of dif\/ferential equation of f\/inite type is based on the general theory of geometric structures on f\/iltered manifolds. First it was developed by Tanaka in \cite{tan70, tan79}. Recall that \emph{a filtered manifold} is a smooth manifold~$M$ equipped with a f\/iltration of the tangent bundle~$TM$ compatible with the Lie bracket of vector f\/ields. At any point $x\in M$ the associated graded vector space $\gr T_xM$ can be endowed with a Lie algebra structure. This nilpotent Lie algebra $\m$ is called \emph{a symbol of a~filtered manifold} (at the point~$x$).  In the paper we consider only the so-called f\/iltered mani\-folds of \emph{constant type}, assuming that the graded nilpotent Lie algebras $\gr T_xM$ are isomorphic to each other for all points $x\in M$.

By \emph{a symbol of a geometric structure on~$M$} we understand a graded Lie algebra $\g$ with the negative part $\g_{-}=\sum_{i<0} \g_i$ which is equal to the symbol~$\m$ of the f\/iltered manifold $M$ of constant type.
Here the Lie algebra $\mathfrak{m}$ is the subalgebra of a so-called universal Tanaka prolongation $\g(\m)$. Roughly speaking, this means that $\g(\m)$ is the maximum among graded Lie algebras which satisfy the condition
``for any element $X\in \g_i$, $i\ge 0$ the equality $[X,\g_-]=0$ implies $X=0$''.

 An arbitrary equation $\E$ can be viewed as a~surface in jet space. The canonical restriction of the contact distribution on jet space def\/ines the structure of f\/iltered manifold on $\E$.

  \subsection{The problem of equivalence} One of the main problems in the theory of dif\/ferential equations is the problem of equivalence.
Two dif\/ferential equations are called equivalent  if one can be transformed to another by a certain change of variables. We consider equations up to point transformations, i.e.\ we allow arbitrary changes of both dependent and independent variables.

First classical approach to the equivalence problem of ODEs was developed by Sophus Lie. In~\cite{lie} he obtains partial results about second order ODEs. The complete answer was given later by Tresse~\cite{tresse}. Invariants of the third order ODEs were computed by Chern in his paper~\cite{chern}. A~modern approach to the equivalence problem of ODEs can be found in the papers~\cite{dou} and~\cite{fels}, where characteristic Cartan connection  was constructed for the one equation of arbitrary order  and for the system of ODEs of the second order.

The general approach to the equivalence problem for the holonomic dif\/ferential equations can be f\/ind in~\cite{dkm}. The key fact there is the existence of a full functor from the category
of holonomic dif\/ferential equations to the category of Cartan connections. This reduces the equivalence problem for dif\/ferential equations to the equivalence problem for the corresponding Cartan connections.

\subsection{Normalization of Cartan connections}

Let $P$ be the principal $H$-bundle. Let $\omega$ be a Cartan connection of type $(G,H)$, where $G$ is a Lie group with a semisimple graded Lie algebra $\g$ and $H$ is a parabolic subgroup of $G$ with the Lie algebra $\gh$. In the paper~\cite{tan79} Tanaka built a set of normal Cartan connections on the principal bundle $P$ as follows. He used the scalar product def\/ined with the help of the Killing form to construct adjoint Lie algebra codif\/ferential $\partial^*$. Then a Cartan connection is normal if\/f the structure function $C:P\to\Hom(\wedge^2 \g_-,\g)$ belongs to the kernel of the operator $\partial^*$ and the structure function has not negative components. As usual def\/ine a Laplacian $\Delta=\partial^*\partial+\partial\partial^* $. The structure function $C$ decomposes as $C=H(C)+\Delta (C)$. The component $H(C)$ is called the harmonic part of the structure function. The key fact about it is that $H(C)$ is the fundamental system of invariants (see Def\/inition~\ref{d2} for details). In the case of the geometry of holonomic dif\/ferential equations the Lie algebra $\g$ is not necessarily semisimple. However in~\cite{dkm} is shown that we still can find the scalar product on $\g$ such that the normal Tanaka conditions def\/ine the unique Cartan connection associated to a holonomic dif\/ferential equation.

In this paper we associate with every system of ODEs of third order a characteristic Cartan connection which dif\/fer from a normal Tanaka Cartan connection. The reason for doing this is a relation between conformal geometry and geometry of the system of the third order ODEs. Conformal manifold is determined by the family of conformal circles, which was shown by Yano~\cite{yano}. Each conformal circle is determined by the point on it, the direction and the curvature, i.e.\ by the point in the third jet space. The system of appropriate dif\/ferential equations of the third order gives us the bridge between the conformal geometry and the geometry of the dif\/ferential equation. It is appeared that a characteristic Cartan connection, which is built in the paper, is in close relations with the normal conformal Cartan connection. The relation of the conformal geometry and the geometry of third order ODEs is the topic of the next paper.

 The paper is organized as follows. In Section~\ref{section2} we naturally associate the system of the third order ODEs with the pair of distributions. This pair of distributions give rise to the f\/iltered manifold associated with the system of the third order ODEs. We write down the symbol of the system of ODEs of the third order, the notion of adopted coframe and adopted Cartan connection.
 The problem of equivalence is considered in Section~\ref{section3}.
 When we working in the case of semisimple Lie algebras and normal Cartan connections, the harmonic part of the curvature gives us the fundamental system of dif\/ferential invariants. We show that in general case fundamental dif\/ferential invariants are contained in the $\ker \partial$ part of the curvature, where $\partial$ is the Lie algebra cohomology dif\/ferential. In Section~\ref{section4} we build the characteristic Cartan connection uniquely associated to the the system of ODEs of the third order.  This connection allows us to obtain the results about equivalence third order equations and to describe the structure of the fundamental invariants of the system of third order ODEs. In particular, this answers the question ``When is the given system trivializable?'' explicitly.

\section{Geometry of the systems of third order ODEs}\label{section2}

Consider an arbitrary system of $m$ ordinary dif\/ferential equations of third order:
 \begin{gather} \label{1}
 y_i'''(x)=f_i\big(y_j''(x),y_k'(x),y_l(x),x\big),
\end{gather}
where $ i,j,k,l=1,\dots,m$ and $m\ge2 $.

We associate a f\/iltered manifold with this system in the following way. Let $J^3(\R^{m+1},1)$ be the third jet space of unparametrized curves. Then the equations~\eqref{1} can be considered as a~submanifold $\E$ in $J^3(\R^{m+1},1)$. We introduce the following coordinate system on the surface~$\E$:
\[(x,y_1,\dots,y_m,p_1=y_1',\dots,p_m=y_m',
q_1=y_1'',\dots,q_m=y_m''). \]

There is a natural one-dimensional distribution $E$ whose integral curves are the lifts of solutions of equations~\eqref{1}. Let $\pi^2_1$ be the canonical projection from the surface $\E$ to the f\/irst jet space
$J^1(\R^{m+1},1)$.  We denote a kernel of a dif\/ferential $\dd\pi^2_1$ as $V$.  In coordinates distributions~$E$,~$V$ have the form:
\begin{gather*}
E=\left\langle \frac{\partial}{\partial x}+p_i\frac{\partial}{\partial y_i}+ q_i\frac{\partial}{\partial p_i} + f^i\frac{\partial}{\partial q_i}\right\rangle,\qquad
V=\left\langle \frac{\partial}{\partial q_i}\right\rangle,
\end{gather*}
where $i,j=1,\dots,m$.

Def\/ine a distribution $C$ as the direct sum of the distributions $E$ and $V$. Then $C$ and its subsequent brackets def\/ine a f\/iltration of a tangent bundle $T\E$:
  \[C=C^{-1}\subset C^{-2} \subset C^{-3} = T\E ,
  \]
where $C^{-i-1}=C^{-i}+[C^{-i},C^{-1}]$.

It is easy to see that the symbol of the f\/iltrated manifold $\E$ is a nilpotent Lie algebra $\m$  isomorphic to the Lie algebra of vector f\/ields generated by
\[
\m_{-1}=\left\langle \frac{\partial}{\partial x}+p_j\frac{\partial}{\partial y_j}+ q_j\frac{\partial}{\partial p_j} , \frac{\partial}{\partial q_i}\right\rangle.
\]

Let $\operatorname{Aut}_0(\m)$ be a subgroup of grading preserving elements of the group $\operatorname{Aut}(\m).$
The elements of the group $\operatorname{Aut}_0(\m)$ which preserve the splitting $ E \oplus V $ form subgroup $G_0$. So the splitting $ E \oplus V $ of the distribution $C$ def\/ines $G_0$-structure of type~$\m$. The action of the group $G_0$ on $\m$ is completely determined by its action on $\m_{-1}$. The latter has the following form in the basis $\left\{ \frac{\partial}{\partial x}+p_j\frac{\partial}{\partial y_j}+ q_j\frac{\partial}{\partial p_j} , \frac{\partial}{\partial q_i}\right\}  $:
\[
    \begin{pmatrix}
      a & 0 \\
      0 & B
    \end{pmatrix} , \qquad a\in\R^* , \quad B\in GL_m(\R). \]

   The symbol $ \g $ is the universal Tanaka prolongation of the pair $(\m,\g_0)$. It has the following form:
  \begin{gather*}
\g=(\gsl _2(\R)\times \gl_m(\R)) \rightthreetimes (V_2 \otimes W).
\end{gather*}
In other words,  $\g$  is equal to the semidirect product of the Lie algebra $\gsl_2(\R)\times
\gl_m(\R)$ and an Abelian ideal $V$. The ideal $V$ has the form $V_2\otimes W$,
where $V_2$ is an irreducible $\gsl_2$-module of dimension $3$ and $W=\R^m$ is the standard representation of $ \gl_m(\R)$.

Let us f\/ix a basis of the Lie algebra $\gsl_2$ and $\gsl_2$-module $V_2$. Let $x$, $y$, $h$ be the
standard basis of an algebra $\gsl_2$ with relations:
\begin{gather*}
  [x,y]=h, \qquad [h,x]=2x,\qquad [h,y]=-2y.
\end{gather*}
This basis can be represented in the following way:
\begin{gather*}
x=
    \begin{pmatrix}
      0 & 1 \\
      0 & 0
    \end{pmatrix},\qquad
  h=
    \begin{pmatrix}
      1 & 0 \\
      0 & -1
    \end{pmatrix},\qquad
y=
    \begin{pmatrix}
      0 & 0 \\
      1 & 0
    \end{pmatrix}.
\end{gather*}
Let $v_0$, $v_1$, $v_2$ be a basis of the module $V_2$ such that $x. v_2=v_1$, $x.v_1=v_0$, $x.v_0=0$.

Def\/ine the grading of the Lie algebra $\g$ as follows:
 \begin{gather*}
\g_1 =\langle y \rangle, \qquad \g_0=\langle h,\gl_m \rangle,\qquad \g_{-1}=\langle
x
\rangle
+\langle v_2\otimes W \rangle,\\
 \g_{-2}=\langle v_1\otimes W \rangle, \qquad
 \g_{-3}=\langle v_0\otimes W \rangle.
\end{gather*}

To build a natural Cartan geometry associated to the equation \eqref{1} we will use the fact~\cite{mor93} that under some additional conditions (which are satisf\/ied for geometric structures arising from holonomic dif\/ferential equations, see~\cite{dkm}) there exists a full functor from the category of $G_0$-structures of type $\m$ to the category of Cartan connections of type $(G,H)$, where $G$ and $H$ are the Lie
groups with Lie algebras $\g$ and $\gh$ respectively which are determined from~$G_0$ in
natural manner. The group $G$ is a semisimple product:
\[
G=\left( SL_2(\R)\times GL_m(\R) \right) \rightthreetimes \left( V_2 \otimes
W \right) , \]
 and the group $H$ is the following subgroup of $G$:
\[ H=  \begin{pmatrix} a & b \\ 0 & a^{-1} \end{pmatrix} \times A, \qquad
a \in \R ^*, \quad b \in \R , \quad A \in GL_m(\R).\]
Note that the corresponding subalgebra $\gh$ is exactly the nonnegative part of the Lie algebra $\g$: \[\gh=\sum_{i\ge 0}\g_i.\]

\begin{definition}
We say that a coframe $\{\omega_{-3}^i,\omega_{-2}^i,\omega_{-1}^i,\omega_x \}$ on $\E$ \emph{is adapted to equation}~\eqref{1} if:
\begin{itemize}\itemsep=0pt
 \item the annihilator of forms $\omega_{-3}^i$, $\omega_{-2}^i$, $\omega_x$ is $V$;
 \item the annihilator of forms $\omega_{-3}^i$, $\omega_{-2}^i$, $\omega_{-1}^i$ is $E$;
 \item the annihilator of forms $\omega_{-3}^i$ is $C^{-2}$.
\end{itemize}
\end{definition}

Let $\overline{\pi}\colon P\rightarrow\E$ be a principle $H$-bundle and let $\overline\omega$ be and arbitrary Cartan connection of type $(G,H)$ on~$P$. Connection $\overline \omega$ can be written as:
\[
 \overline{\omega }=\overline{\omega }^i_{-3} v_0\otimes e_i +
 \overline{\omega }^i_{-2} v_1\otimes e_i +
 \overline{\omega }^i_{-1} v_2\otimes e_i + \overline{\omega }_x x + \overline{\omega }_h h +
\overline{\omega }^i_j e^j_i + \overline{\omega }_y y .
\]

\begin{definition}
 We say that a Cartan
connection $\overline{\omega}$ on a principal $H$-bundle $\overline{\pi}$ is adapted to equations~\eqref{1}, if
for any local section $s$ of $\overline{\pi}$ the set
 \[
 \big\{ s^* \overline{\omega}_x , s^*\overline{\omega}_{-1}^i, s^*\overline{\omega}_{-2}^i,
s^*\overline{\omega}_{-3}^i \big\}
\] is an adapted co-frame on~$\E$.
\end{definition}

We have described the set of Cartan connection adapted to the system of third order ODEs. However, we can chose the representative in dif\/ferent ways. The next two sections are devoted to the building of a~canonical connection which we call characteristic.

\section{Characteristic Cartan connection\\ and fundamental dif\/ferential invariants}\label{section3}

As in Section~\ref{section2} let $\overline{\pi}\colon P\rightarrow\E$ be a principle $H$-bundle and let $\overline\omega$ be and arbitrary Cartan connection of type $(G,H)$ on $P$:
\[
 \overline{\omega }=\overline{\omega }^i_{-3} v_0\otimes e_i +
 \overline{\omega }^i_{-2} v_1\otimes e_i +
 \overline{\omega }^i_{-1} v_2\otimes e_i + \overline{\omega }_x x + \overline{\omega }_h h +
\overline{\omega }^i_j e^j_i + \overline{\omega }_y y .
\]
Let $\overline{\Omega}=\dd \overline{\omega}+\frac{1}{2}[\overline{\omega},\overline{\omega}]$ be the
 curvature of the Cartan connection $\overline\omega$:
\[ \overline{\Omega}=\overline\Omega_{-3}^i v_0\otimes e_i+ \overline\Omega_{-2}^i v_1\otimes e_i+
\overline\Omega_{-1}^i
v_2\otimes e_i+  \overline\Omega_x x+ \overline\Omega_h h+ \overline\Omega_j^i e_i^j+\overline\Omega_y y .\]
\begin{definition}
The structure function of a Cartan connection $\omega$ is  a function
\[C:P\to \Hom\big({\wedge}^2 \g_-, \g\big),\]
which is def\/ined by
\[C(p)( g_1, g_2) = \Omega_p \big(\omega_p^{-1}(g_1) , \omega_p^{-1}(g_2) \big).\]
\end{definition}

We can obtain the structure function of a Cartan connection explicitly. Let $\{e_1,\dots, e_{n+k} \}$ be a basis of  Lie algebra $\g$ such that
$\{e_{n+1} , \dots , e_{n+k} \}$ form a basis of the subalgebra $\gh.$ In our case $\{e_{n+1} , \dots , e_{n+k} \}=\{h,y,e^j_i\}.$ An arbitrary element $\varphi\in\Hom(\wedge^2 \g_-, \g)$ def\/ined by constants~$C^k_{ij}$, where
\[ \varphi(e_i,e_j)=\sum_{k=1}^{n+k}C^k_{ij}e_k,\qquad 1 \le i, j \le n.
 \]
 The structure function $C : P \to \Hom(\wedge^2 \g_-,\g)$ def\/ines  functions $ C^k_{ij} (p).$ If
 \[\omega=\sum \omega_i e_i,\qquad\Omega=\sum \Omega^k e_k,\]
 then the functions $ C^k_{ij} (p)$ can be found from the decomposition of the curvature tensor $\Omega $ in terms of forms $\omega_i$:
 \[ \Omega^k = \sum C^k_{ij} \omega_i\wedge\omega_j. \]

Let $\Omega^i$ be one of the 2-forms $\overline\Omega_{-3}^i$, $\overline\Omega_{-2}^i$,
$\overline\Omega_{-1}^i$,
$\overline\Omega_x$, $\overline\Omega_h$, $\overline\Omega_j^i$.
We can write it explicitly as:
\[
\Omega^i=\sum_{p,q=1}^3\Omega^i\big[\overline\omega_{-q}^j,\overline\omega_{-p}^k\big]\overline\omega_{-q}^j
\wedge \overline{\omega}_{-p}^k+ \sum_{p=1}^3 \Omega^i\big[\overline\omega_x,\overline\omega_{-p}^k\big]\overline\omega_x \wedge \overline{\omega}_{-p}^k.
\]
Then $\Omega^i[\overline\omega_{-q}^j,\overline\omega_{-p}^k]$ and $ \Omega^i[\overline\omega_x,\overline\omega_{-p}^k]$ are the coef\/f\/icients of the structure function of the Cartan connection $\omega.$ The grading of Lie algebra $\g$ induces degree of the coef\/f\/icients $\Omega^i[\overline\omega_{-q}^j,\overline\omega_{-p}^k]$ and $ \Omega^i[\overline\omega_x,\overline\omega_{-p}^k]$.

\begin{definition}
We say that Cartan connection associated with the equation~\eqref{1} is characteristic if the following conditions on a  curvature is satisf\/ied:
\begin{itemize}\itemsep=0pt
 \item all coef\/f\/icients  of degree $\le 1$ are equal to $0$;
 \item in degree $2$ we have $\overline\Omega_h[\overline\omega_x\wedge \overline\omega_{-1}^i]=0$,  $\overline\Omega^i_j[\overline\omega_x\wedge
\overline\omega_{-1}^k]=0$,  $\overline\Omega_x[\overline\omega_x\wedge \overline\omega_{-2}^i]=0$, $\overline\Omega_{-1}^i [\overline\omega_x\wedge
\overline\omega_{-2}^i]=0$;
 \item in degree $3$ we have $\overline\Omega_y[\overline\omega_x\wedge \overline\omega_{-1}^i]=0$, $\overline\Omega_h[\overline\omega_x\wedge
\overline\omega_{-2}^i]=0$, $\overline\Omega^i_j[\overline\omega_x\wedge\overline\omega_{-2}^k]=0$;
\item in degree $4$ we have $\overline\Omega_y[\overline\omega_x\wedge \overline\omega_{-2}^i]=0$.
\end{itemize}
\end{definition}

In other worlds these conditions def\/ine the subspace $U$ and Cartan connection is characteristic if and only if it belongs to $U$.
\begin{theorem}\label{t1}
 There exists a unique characteristic Cartan connection
associated with the equation~\eqref{1}.
\end{theorem}

\begin{proof}
We will proceed with parametric computations of characteristic Cartan connection in the forth section of the paper. We will f\/ix a section $s:\E \to P$ and prove that locally for every equation there exists a unique Cartan connection $\omega$ with structure function pullback $ s^*C:\E\to\Hom(\wedge^2 \g_-,\g)$ takes values in the space $U$. Now we show that the characteristic Cartan connection is uniquely globally def\/ined with this data.

Take a covering $U_\alpha $ of the space $E$
and construct a Cartan connection $\omega_\alpha$ on each trivial f\/ibre bundle $\pi_\alpha : U_\alpha \times H \to U_\alpha$. Let $s_\alpha$ and $s_\beta$ be the trivial sections of the f\/ibre bundles~$\pi_\alpha$ and~$\pi_\beta$. Let $\widetilde\omega_\alpha=s_\alpha^*\omega_\alpha$ and $\widetilde\omega_\beta=s_\beta^*\omega_\beta$. Since forms $\omega_\alpha$ and $\omega_\beta$ are uniquely def\/ined there exists a~unique function
\[\varphi_{\alpha\beta} : U_\alpha \cap U_\beta \to H ,\]
such that
\[\omega_\beta = \operatorname{Ad}\big( \varphi^{-1}_{\alpha\beta} \big) 
\omega_\alpha+ \varphi^*_{\alpha\beta}\omega_H ,\]
where $ \omega_H$ is Maurer--Cartan form of the Lie group $H$. The functions $\varphi_{\alpha\beta}$ uniquely def\/ine a~principle $H$-bundle with the Cartan connection $\omega$.

In order to prove that the structure function $C$  of the Cartan connection $\omega$ takes values in the space~$U$ it is suf\/f\/icient to show that~$U$ is $\operatorname{Ad}(H)$-invariant.

Note that the action of $G_0$ preserves the zero condition on the structure function of the characteristic connection. We need only to check that the space $U$ is $\operatorname{exp}(y)$-invariant or equally $\operatorname{ad}(y)$ invariant. The action of the element $y$ has degree one. Conditions on the curvature of the Theorem \ref{t1} are $\operatorname{ad}(y)$-invariant up degree 2, since all components of degree less than~2 are equal to zero. Finally, the conditions of degree 3 and 4 are $\operatorname{ad}(y)$-invariant, since the coef\/f\/icients $\overline\Omega_y[\overline\omega_x\wedge \overline\omega_{-1}^i]$, $\overline\Omega_h[\overline\omega_x\wedge
\overline\omega_{-2}^i]$, $\overline\Omega^i_j[\overline\omega_x\wedge\overline\omega_{-2}^k]$ and $\overline\Omega_y[\overline\omega_x\wedge \overline\omega_{-2}^i]$ can be obtained only from
 $\overline\Omega_h[\overline\omega_x\wedge \overline\omega_{-1}^i]$, $\overline\Omega_x[\overline\omega_x\wedge
\overline\omega_{-2}^i]$,  $\overline\Omega^i_j[\overline\omega_x\wedge\overline\omega_{-1}^k]$, $\overline\Omega_h[\overline\omega_x\wedge \overline\omega_{-2}^i]$ and  $\overline\Omega_y[\overline\omega_x\wedge \overline\omega_{-1}^i]$ which all are zero for characteristic Cartan connection. This ends the proof of a global existence of the form $\omega$.
\end{proof}

{\sloppy Let $V$ be an arbitrary f\/inite-dimensional vector space and let $f$ be a smooth function \mbox{$f:P\to V$}. Denote by $L_0(f)$ the space of all functions of the form $\langle f, v^*\rangle $,
where $v^* \in V^*$ and by~$L(f)$ the algebra generated by elements from $L_0 (f )$ and all their covariant derivatives. For example, the algebra $L(C)$, where $C$ is structure function of the Cartan connection $\omega$, consists of local invariants of the connection~$\omega$.

}

\begin{definition}\label{d2}
We say that functions $f_i$ are the fundamental system of dif\/ferential invariants for the structure with Cartan connection $\omega$ if $L(f_i)=L(C)$.
\end{definition}

 The key to calculation of the fundamental system of dif\/ferential invariants is to determine which parts of the curvature are expressed through another. In~\cite{dkm} it is shown that fundamental invariants of holonomic dif\/ferential equation lie in non-negative harmonic part of the curvature of the normal Cartan connection. In general we have approximately the same situation: there is one to one correspondence between fundamental dif\/ferential invariants of the characteristic Cartan connection and $H^2_+(\g_-,\g)$ part of the structure function. Here $H^2_+(\g_-,\g)$ is the non-negative part of the second Lie algebra cohomology group.

\begin{proposition}\label{p4}
Let $\omega$ be a Cartan connection of type $(G,H)$ on a principal $H$-bundle $P$, where $(G,H)$ is an arbitrary pair of Lie group and its subgroup. Assume that the  Lie algebra $\g$ is a~graded Lie algebra of the Lie group $G$ with the negative part~$\g_-$. Assume that a structure function of $\omega$ takes values in subspace $W\subset\Hom(\wedge^2\g_-,\g)  $ and has only components of positive degree. Then a   $\ker \partial \cap W$ part of the structure function forms a system of fundamental differential invariants.
\end{proposition}

\begin{proof}
The algebra of dif\/ferential invariants is generated by the structure function coef\/f\/icients. We will use the Bianchi identity to show that some coef\/f\/icients of the characteristic Cartan connection curvature are obtained from the image of the operator~$\partial$.

 Let $e_i$ be the basis of the Lie algebra $\g$, $X_i$ be the corresponding fundamental vector f\/ields on $P$ and $\omega^i$ be the dual coframe. We can write  the Cartan connection $\omega$ in the form:
\[ \omega= \omega^i e_i.\]
Assume that the Lie algebra $\g$ has structure constants $A_{ij}^k$. That means that:
\[ [e_i,e_j]=A_{ij}^k e_k .\]
Write the curvature of the Cartan connection $\omega$ in coordinates:
\begin{gather}\label{41}
 \Omega=C_{ij}^k \omega^i \wedge \omega^j e_k.
 \end{gather}
Then the following equality is fulf\/illed:
\[ \dd \omega^k =\big(C^k_{ij}- A_{ij}^k\big) \omega^i \wedge \omega^j .\]
Now apply the Bianchi identity $\dd \Omega=[\Omega,\omega]$ to the equation \eqref{41}:
\[ \left(\frac{\partial C_{ij}^k}{\partial X_l}\omega_l\wedge\omega_i\wedge\omega_j e_k+C_{ij}^k \dd \omega_i\wedge\omega_j e_k+ C_{ij}^k \omega_i\wedge \dd\omega_j e_k\right)=
C_{ij}^k[e_k,e_l]\omega_i\wedge\omega_j\wedge\omega_l.
\]
Express the covariant derivative of the structure function:
\begin{gather*}
\frac{\partial C_{ij}^p}{\partial X_l}\omega_l\wedge\omega_i\wedge\omega_j e_p=
 \big({-}C_{kl}^p \dd \omega_k\wedge\omega_l-C_{kl}^p \omega_k\wedge \dd\omega_l +C_{ij}^k A^p_{kl}\omega_i\wedge\omega_j\wedge\omega_l\big)e_p \\
\hphantom{\frac{\partial C_{ij}^p}{\partial X_l}\omega_l\wedge\omega_i\wedge\omega_j e_p}{}  =
 C^p_{kl}\big(C^k_{ij}-A^k_{ij}\big) \omega_i\wedge\omega_j\wedge\omega_l e_p +  C_{ij}^k A^p_{kl}\omega_i\wedge\omega_j\wedge\omega_l e_p
.
\end{gather*}
We get that:
\begin{gather}
 \frac{\partial C_{ij}^p}{\partial X_l}\omega_l\wedge\omega_i\wedge\omega_j e_p- C^p_{kl}C^k_{ij} \omega_i\wedge\omega_j\wedge\omega_l e_p\nonumber\\
 \qquad{} =
  C^p_{kl}A^k_{ij} \omega_i\wedge\omega_j\wedge\omega_l e_p + C_{ij}^k A^p_{kl}\omega_i\wedge\omega_j\wedge\omega_l e_p .\label{42}
\end{gather}
If we take  the $\Hom(\wedge^3 \g_-,\g)$ part of \eqref{42} (i.e.\ assume that $\omega_l\in\g_-^*$) we get that the right side of the \eqref{42} is exactly the Lie cohomology dif\/ferential.

On the right side of \eqref{42} coef\/f\/icients have the same degree as in the curvature. On the other hand  coef\/f\/icients on the left side have an increased degree. So, we have obtained that coef\/f\/icients which are mapped to the $\operatorname{im} \partial$ can be expressed through the covariant derivative of the coef\/f\/icients of the lower degree. This proves the proposition.
\end{proof}

\begin{remark}
Note that if intersection of  $W$ and $\im \partial $ is zero then subspace $\ker \partial\cap W$ is generated by  representatives of  $H^2_+(\g_-,\g)$.
\end{remark}

\begin{theorem}\label{t2}
The following invariants are fundamental differential invariants for the system of third order ODEs:
\begin{gather*}
\left(W_2\right)^i_j  = \tr_0 \left(\frac{\partial f^i}{\partial p^j} -
\frac{\dd}{\dx} \frac{\partial f^i}{\partial q^j} +
\frac{1}{3} \frac{\partial f^i}{\partial q^k}\frac{\partial f^k}{\partial q^j} \right), \qquad
\left(I_2\right)^i_{j,k}  = \tr_0\left(\frac{\partial^2 f^i}{\partial q^j \partial
q^k}
 \right), \\
\left(W_3\right)^i_j  =  \frac{\partial f^i}{\partial y^j} +\frac{1}{3} \frac{\partial f^i}{\partial q^k}\frac{\partial f^k}{\partial p^j} -\frac{1}{2} \frac{d}{dx} \frac{\partial f^i}{\partial
p^j} + \frac{1}{6}\frac{d^2}{dx^2} \frac{\partial f^i}{\partial
q^j} - \frac{2}{27}\left(\frac{\partial f^i}{\partial q^k}\right)^3- \frac{1}{18}
\frac{\partial f^i}{\partial q^k}\frac{d}{dx}\frac{\partial f^k}{\partial q^j}
\\ \phantom{\left(W_3\right)^i_j}{} -
\frac{5}{18} \frac{d}{dx}\left(\frac{\partial f^i}{\partial q^k}\right)\frac{\partial f^k}{\partial q^j} ,\\
 \left(I_4\right)_{j,k}  = - \frac{\partial H_k^{-1}}{\partial p_j}+\frac{\partial}{\partial q_j}\frac{\partial}{\partial q_k}H^x - \frac{\partial}{\partial q_k}\frac{\dd}{\dx}H_j^{-1} - \frac{\partial }{\partial q^k}\left(H_l^{-1}\frac{\partial
f^l}{\partial q^j}\right)+2H_j^{-1} H_k^{-1},
\end{gather*}
where
\[
H_j^{-1}=\frac{1}{6(m+1)}\left(\frac{\partial^2 f^i}{\partial q^i \partial
q^j}
 \right) \qquad \mbox{\rm and} \qquad H^x=-\frac{1}{4m}\left(\frac{\partial f^i}{\partial p^i} -
\frac{\dd}{\dx} \frac{\partial f^i}{\partial q^i} +
\frac{1}{3} \frac{\partial f^i}{\partial q^k}\frac{\partial f^k}{\partial q^i} \right).
\]
\end{theorem}

\begin{proof}
We will use Proposition~\ref{p4}. The fundamental dif\/ferential invariants is in one to one correspondence with the cohomology group $H^2_+(\g_-,\g)$.  For the case of the system of ODEs of the third order the Lie  cohomology group $H^2_+(\g_-,\g)$ was studied in~\cite{medv}.  The main result of that work is that the space $H^2_+(\g_-,\g)$ has the following decomposition:
\begin{center}
 \begin{tabular}{|c|c|}
 \hline  Degree & Space  \\
\hline
$-1$ & $v^0_6 \otimes \wedge^2(W^*) \otimes W$\tsep{2pt}\\
$0$ & $v^0_4 \otimes S_0^2(W^*) \otimes W$  \\
$0$ & $v^0_4 \otimes \wedge^2(W^*) \otimes W$ \\
$1$ & \small $v^0_2 \otimes \wedge^2W^*\otimes W/V_2 \otimes W^*$ \\
$2$ & \small $x^* \otimes \R y \otimes \gsl(W) $ \\
$2$ & \small $v_0^0 \otimes S^2(W^*) \otimes W$  \\
$3$ & \small $ x^* \otimes \R y^2 \otimes \gl(W)$ \\
$4$ & \small $v_0^0 \otimes S^2(W^*)$  \\
$3$ & \small $v_2^0 $ if $ m=2 $  \\
\hline
\end{tabular}
\end{center}
Here $v_k^0$ is the lowest vector of corresponding $(k+1)$-dimensional $\gsl_2$-module $V_k$.

 Now we list the result table with the corresponding invariant. We start from degree 2 since all part of curvature of degree less than 2 is zero.
\begin{center}
 \hspace*{25mm}\begin{tabular}{|c|c|c|c|}
 \hline  Degree & Space & Part of the curvature & Invariant \\
\hline
$2$ & \small $x^* \otimes \R y \otimes \gsl(W)$ &$\Omega_{-1}^i [\omega_x\wedge\omega_{-2}^j]$ & $W_2$\tsep{2pt}
 \\
$2$  & \small $v_0^0 \otimes S^2(W^*) \otimes W$ &$\Omega_{-2}^i [\omega_{-1}^j\wedge\omega_{-3}^k]$ &$I_2$ \\
$3$ & \small $ x^* \otimes \R y^2 \otimes \gl(W)$ &$\Omega_{-1}^i [\omega_x\wedge\omega_{-3}^j]$&$W_3$ \\
$4$  & \small $v_0^0 \otimes S^2(W^*)$  &$\Omega_{y} [\omega_x\wedge\omega_{-3}^j]$ &$I_4$ \\
$3$ & \small $v_2^0 $  if $ m=2 $ & $\Omega_{y} [\omega_{-1}^2\wedge\omega_{-2}^1]$&$\equiv 0$ \\
\hline
\end{tabular}\hfill\qed
\end{center}
\renewcommand{\qed}{}
\end{proof}

\begin{corollary}
The system \eqref{1} is equivalent to the trivial one via point transformations if and only if all invariants $I_2$, $W_2$, $W_3$, $I_4$ vanish identically.
\end{corollary}

\begin{example}[Dif\/ferential equations on circles in $\R^n$]
As application of the previous results we compute invariants of the system of third order ODEs on circles in Euclidean space.

\begin{lemma}
 Let  $ E$ be the $(m+1)$-dimensional Euclidean space with the orthonormal basis
 $\{ e_0, \ldots, e_n\}$ and the coordinates $\{r_0,r_1,\dots,r_n\}$. Then the equation of circles in $E$ parametrized by the coordinate $r_0$ is:{\samepage
\[\dddot{r}_i = 3\ddot{r}_i \frac{\sum\limits_{j=1}^m \dot{r}_j \ddot{r}_j } {1+ \sum\limits_{j=1}^m \dot{r}^2_j} , \qquad i=1,\ldots,m  .\]
This equation is invariant under conformal transformations of~$E$.}
\end{lemma}

\begin{proof}
Let the curve $R(t)= (r_0(t), \dots, r_n(t))$ be a circle. Assume now that $r_0(t)=t$. We have
\begin{gather}
  \dddot{R}(t) = a(t)\ddot{R}(t) + b(t)\dot{R}(t), \label{102}
\end{gather}
since $R(t)$ is 2-dimensional curve.
Next, $ b(t) = 0 $ in our parametrization, since
\[ 0=\dddot {r}_0(t) = a(t)\ddot{r}_0(t) + b(t) \dot{r}_0(t) =b(t).\]
To determine $a(t)$ note that
 \[(R(t)-C,R(t)-C)=d \]
 for some constant $d$ and $C\in E$. Dif\/ferentiating, we get:
 \begin{gather*}
(\dot{R}(t),R(t)-C)=0 , \\
(\ddot{R}(t),R(t)-C)=-(\dot{R}(t),\dot{R(t)}), \\
(\dddot{R}(t),R(t)-C)+3(\ddot{R}(t),\dot{R}(t))=0.
\end{gather*}
Now substitute \eqref{102} into previous formula:
\begin{gather*}
(a(t)\ddot{R}(t)+b(t)\dot{R}(t), R(t)- C)=-3(\ddot{R}(t),\dot{R}(t)), \\
 (a(T)\ddot{R}(t), R(t)-C)=-a(t)(\dot{R}(t),\dot{R}(t))=-3(\ddot{R}(t),\dot{R}(t)).
 \end{gather*}
We get that \[
 a(t)=3\frac{(\ddot{R}(t),\dot{R}(t))}{(\dot{R}(t),\dot{R}(t))} .
\]
Substituting $ a(t) $ into \eqref{102} we get our equations.
\end{proof}

\begin{proposition}
For differential equation on conformal circles invariants $W_2$, $I_2$, $W_3$ vanish identically. Invariant $I_4$ has the following form:
\[
\left(I_4\right)^i_j =\frac{1}{2} \delta^i_j \frac{1}{1+\sum\limits_{k=1}^m \dot{r}_k^2}- \frac{1}{2}\frac{\dot{r}_i \dot{r}_j}{\left(1+\sum\limits_{k=1}^m \dot{r}_k^2\right)^2}.
\]
\end{proposition}

\begin{proof}
The proof is straightforward applying of the formulas from Theorem \ref{t2}.
\end{proof}
\begin{remark} There are other equations satisfying $W_2=I_2=W_3=0$. For example, it is an union of a system on circles in $R^{n-k}$ and a system of $k$ trivial equations. It would be interesting to characterize geometrically the class of such equations.
\end{remark}
\end{example}

\section{Parametric computation\\ of the characteristic Cartan connection}\label{section4}

Consider a system of third-order ordinary dif\/ferential equations of the form
\begin{gather*}
(y^i)'''=f^i\big(x,y^j,(y^k)',(y^l)''\big) ,
\end{gather*}
where $i,j=1,\dots,m$ with $m\ge2$. It determines a holonomic
dif\/ferential equation $\E\subset J^3(\R^{m+1},1)$. Let us use the following coordinate system on the equation~$\E$:
\[x,y_1,\dots,y_m,p_1=y_1',\dots,p_m=y_m',
q_1=y_1'',\dots,q_m=y_m''. \]
 We choose a coframe $\theta$ on the surface $\E$:
\begin{gather*}
  \theta_x   =dx;\\
  \theta^i_{-1}   =dq^i-f^i(x,y,p,q)\,dx,\quad i=1,\dots,m;\\
  \theta^i_{-2}   =dp^i-q^i\,dx,\quad i=1,\dots,m;\\
  \theta^i_{-3}   =dy^i-p^i\,dx,\quad i=1,\dots,m.
\end{gather*}
To connect our computation on the surface $\E$ with the principle bundle $P$ let us use the following uniquely def\/ined section  $s\colon\E\rightarrow P$ with relations:
\begin{gather*}
 s^* \overline{\omega}_{-3}^i  =\theta_{-3}^i ,\\
 s^* \overline{\omega}_h  \equiv 0 \mod \langle \theta_{-3}^i,\theta_{-2}^i,\theta_{-1}^i \rangle ,\\
 s^* \overline{\omega}_x  \equiv -\theta_x \mod \langle \theta_{-3}^i,\theta_{-2}^i,\theta_{-1}^i \rangle.
\end{gather*}

Def\/ine a pullback  $\omega \colon T\E \to \g$ by the formula $ \omega  = s^* \overline{\omega }$.
Let $\overline{\Omega}$ be a curvature tensor of~$\overline{\omega }$, and let
$\Omega=s^*\overline{\Omega}$.
We see that
\begin{gather*}
\Omega  =\Omega^i_{-3} v_0\otimes e_i + \Omega^i_{-2} v_1\otimes e_i +
\Omega^i_{-1} v_2\otimes e_i + \Omega_x x + \Omega_h h + \Omega_i^j e_j^i + \Omega_y y \\
\phantom{\Omega}{} =(d \omega_{-3}^i + \omega_x \wedge \omega_{-2}^i + 2\omega_h \wedge
\omega_{-3}^i + \omega_j^i \wedge \omega_{-3}^j)  v_0\otimes e_i    \\
\phantom{\Omega=}{} +(d \omega_{-2}^i + \omega_x \wedge \omega_{-1}^i +  \omega_j^i \wedge \omega_{-2}^j + 2 \omega_y \wedge
\omega_{-3}^i)  v_1\otimes e_i    \\
\phantom{\Omega=}{}+ (d \omega_{-1}^i  - 2 \omega_h \wedge
\omega_{-1}^i + \omega_j^i \wedge \omega_{-1}^j + 2 \omega_y \wedge
\omega_{-2}^i)  v_2\otimes e_i   \\
\phantom{\Omega=}{}+ (d \omega _x + 2\omega _h \wedge \omega_x)x + (d \omega _h + \omega _x \wedge \omega_y)h     + (d \omega_j^i  + \omega_k^i \wedge \omega_j^k)e_i^j + (d \omega _y -2 \omega _h \wedge \omega_y)y.
\end{gather*}
An arbitrary Cartan connection adapted to equation \eqref{1} has the form:
\begin{gather*}
\omega_{-3}^i  =\theta_{-3}^i ,\\
\omega_{-2}^i  =\alpha^i_j\theta_{-2}^j + A_j^{i}\theta_{-3}^j,\\
\omega_{-1}^i  =\beta^i_j\theta_{-1}^j + B_j^{i}\theta_{-2}^j+C_j^{i}\theta_{-3}^j,\\
\omega_x  =-\theta_x + D_j\theta^j_{-2}+E_j\theta_{-3}^j , \\
\omega_h  = F^{-1}_j\theta^j_{-1}+F^{-2}_j\theta^j_{-2}+ F^{-3}_j\theta^j_{-3},\\
\omega^i_j  =G^{i,x}_j \theta_x +
G^{i,-1}_{jk}\theta^k_{-1}+G^{i,-2}_{jk}\theta^k_{-2}+G^{i,-3}_{jk}\theta^k_{-3},\\
\omega_y  =H^x \theta_x + H^{-1}_j\theta^j_{-1}+H^{-2}_j\theta^j_{-2}+ H^{-3}_j\theta^j_{-3}.
\end{gather*}

In degree 0 of the curvature we have two nonzero components:
\begin{gather*}
  \Omega_{-3}^i \mod \langle \theta_{-2} \wedge \theta_{-2},\theta_{-3} \rangle
=\theta^x \wedge \theta_{-2}^i-\alpha^i_j\theta^x \wedge \theta_{-2}^j, \\
  \Omega_{-2}^i \mod \langle \theta_{-2} ,\theta_{-3} \rangle
=\theta^x \wedge \theta_{-1}^i-\beta^i_j\theta^x \wedge \theta_{-1}^j .
\end{gather*}
Assume these two equalities is zero and get $\alpha^i_j=\delta_j^i$ and $\beta_j^i=\delta_j^i$.

We have three nonzero components in degree 1. The f\/irst component is:
\begin{gather*}
  \Omega_{-3}^i \mod \langle
\theta_{-2} \wedge \theta_{-3},\theta_{-3} \wedge \theta_{-3}  \rangle\\
\qquad {}  =
  -\theta_x \wedge  A^i_j\theta_{-2}^j +D_j \theta_{-2}^j \wedge \theta_{-2}^i+G_i^{i,x}\theta_x \wedge
\theta_{-3}^j+G_{jk}^{i,-1}\theta_{-1}^k\wedge \theta_{-3}^j+2F_j^{-1}\theta_{-1}^j \wedge
\theta_{-3}^i .
\end{gather*}
The second component is:
\begin{gather*}
 \Omega_{-2}^i \mod \langle \theta_{-2} \wedge \theta_{-2},\theta_{-3} \rangle     \\
\qquad{}=A^i_j\theta_x \wedge \theta^j_{-2}+D_j\theta_{-2}^j \wedge \theta^i_{-1}-\theta_x \wedge
B^i_j\theta^j_{-2}+G^{i,x}_j\theta_x \ \wedge \theta^j_{-2}+G^{i,-1}_{jk}\theta^k_{-1} \wedge
\theta^j_{-2}.
\end{gather*}
The third component is:
\begin{gather*}
  \Omega_{-1}^i \mod \langle \theta_{-2} ,\theta_{-3} \rangle  \\
\qquad{}= \frac{\partial f^i}{\partial
q^j}\theta_x \wedge \theta^j_{-1}+B^i_j\theta_x \wedge \theta^j_{-1} -2
F^{-1}_j\theta^j_{-1} \wedge \theta^i_{-1}+ G^{i,x}_j\theta_x \wedge \theta^j_{-1}
+G^{i,-1}_{jk}\theta^k_{-1}\wedge \theta^j_{-1} .
\end{gather*}
After applying zero conditions to these parts of the curvature we obtain
\[
A_j^i=G^{i,x}_j=\frac{1}{2}B^i_j=-\frac{1}{3}\frac{\partial f^i}{\partial q^j},\qquad
D_j=F_j^{-1}=G^{i,-1}_{jk}=0.
\]

Proceed now to the second degree
\begin{gather*}
 \Omega_{-1}^i \mod \langle \theta_{-2}\wedge \theta_{-2},\theta_{-3} \rangle  \\
\qquad{}=\frac{\partial f^i}{\partial p_j}\theta_x \wedge \theta^j_{-2}+ 2\frac{d A^i_j}{d x}\theta_x \wedge
\theta^j_{-2}+  2\frac{\partial  A^i_j}{\partial q_k}\theta^k_{-1} \wedge \theta^j_{-2} +C^i_j \theta_x
\wedge \theta^j_{-2}-2F^{-2}_j \theta^j_{-2}\wedge \theta^i_{-1}
  \\
\qquad\quad{}+G^{i,-2}_{jk}\theta^k_{-2}\wedge  \theta^j_{-1}+2 H^x\theta_x\wedge
\theta^i_{-2}+2H^{-1}_j\theta^j_{-1}\wedge \theta
^i_{-2}+G^{i,x}_k\theta_x \wedge  B^k_j \theta^j_{-2}.
\end{gather*}
We have:
\[\Omega_{-1}^i \big[  \theta _x \wedge \theta^j_{-2}\big]= \frac{\partial f^i}{\partial p_j} +2\frac{d A^i_j}{d
x}+C^i_j+2H^x+   2A^i_k A^k_j .
 \]
 Assuming the previous tensor
is zero, we obtain:
\[
C^i_j=-\left(\frac{\partial f^i}{\partial p_j} +2\frac{d A^i_j}{d
x}+2H^x+   2A^i_k A^k_j \right).
\]
Next curvature component  contains all second order invariants:
 \begin{gather*}
  \Omega_{-2}^i \mod \langle \theta_{-2}\wedge \theta_{-3},\theta_{-3} \wedge \theta_{-3}
\rangle    \\
\qquad{}= \frac{d A^i_j}{d
x}\theta_x \wedge \theta^j_{-3}+ \frac{\partial A^i_j}{\partial q_k}\theta^k_{-1} \wedge \theta^j_{-3}-\theta_x
\wedge C^i_j \theta^j_{-3}+E_j \theta^j_{-3}\wedge \theta^i_{-1}+G^{i,x}_j \theta_x \wedge A^j_k
\theta^k_{-3}  \\
\qquad\quad{}+2H^x\theta_x  \wedge \theta^i_{-3} +2H^{-1}_j\theta^j_{-1}\wedge
\theta^i_{-3}+G^{i,-2}_{jk}\theta^k_{-2} \wedge \theta^j_{-2} +G^{i,-1}_{jk}\theta^{-1}_k \wedge
A^j_l \theta_{-3}^l .
\end{gather*}
In coef\/f\/icient  $\Omega^i_{-2}[\theta ^k_{-1} \wedge \theta^j_{-3}]$ we get invariant $I_2$
\[ \Omega^i_{-2}\big[\theta ^k_{-1} \wedge \theta^j_{-3}\big] = \frac{\partial A^i_j}{\partial q_k}-E_j \delta^i_k +2H^{-1}_k
\delta^i_j= \frac{\partial A^i_j}{\partial q_k} + 2H^{-1}_k \delta^i_j +2F^{-2}_j \delta^i_k .\]
Explicitly, the invariant $I_2$ is the following:
\[I_2=\tr_0 \left( \frac{\partial^2 f^i}{\partial q_j \partial q_k} \right) ,\]
where $\tr_0$ is a traceless part of the tensor.

In the coef\/f\/icient
\[\Omega_{-2}^i\big[\theta_{x} \wedge \theta^j_{-3}\big]=-C^i_j\frac{d A^i_j}{d x}+ A^i_k A^k_j+2H^k \delta^i_j\] we obtain a so-called generalized Wilczynski
invariant.
 As shown in \cite{dou}, a part of dif\/ferential invariants of systems of ODEs comes from its linearisation. As in \cite{dou}, we call them generalized Wilczynski invariants. In our case we have two Wilczynski invariants of degree~2 and~3. We
denote them as $W_2$ and $W_3$ respectively.
The second degree generalized Wilczynski invariant is the following:
 \[  W_2=\tr_0 \left( \frac{\partial f^i}{\partial p_j} -\frac{d }{d x} \frac{\partial f^i}{\partial q_j}
+\frac{1}{3} \frac{\partial f^i}{\partial q_k} \frac{\partial f^k}{\partial q_j} \right).
 \]
Normalizing the trace of previous tensor to zero we obtain:
 \[  H^x=-\frac{1}{4m} \left(\frac{\partial f^i}{\partial p_i}+3\frac{d A^i_i }{d x} +3  A^i_k
A^k_i\right) .\]

 It remains to compute only $\gsl_2\times\gl_m$ part of the curvature in degree 2.
 \[  \Omega_x \mod \langle \theta_{-2}\wedge \theta_{-2},\theta_{-3} \rangle = E_j \theta_x \wedge
\theta^j_{-2}+2 F_j \theta_x \wedge \theta^j_{-2} .\]
Assuming that it vanishes identically we get the following condition:
\[E_j=-2F^{-2}_j.\]
We have:
 \[  \Omega_h \mod \langle \theta_{-2},\theta_{-3} \rangle = F^{-2}_j \theta_x \wedge \theta^j_{-1}- \theta_x
\wedge \theta^j_{-1}H^{-1}_j .\]
The condition $\Omega_{h}^i[\theta _{x} \wedge \theta^i_{-1}]=0$ gives equality $ F^{-2}_j=H^{-1}_j .$

Assuming the trace of the tensor $\Omega_{-2}^i[\theta _{-1}^j \wedge \theta^k_{-3}]$ is equal to zero we get:
  \[
  F^{-2}_k=H^{-1}_k=- \frac{1}{2(m+1)}\frac{\partial A^i_i}{\partial q_k}.
  \]
The last part of degree 2 calculation is:
  \[
  \Omega^i_j \mod \langle \theta_{-2},\theta_{-3} \rangle = \frac{\partial A^i_j}{\partial q_k} \theta_x \wedge
\theta^{-1}_k + G^{i,-2}_{jk}\theta_x \wedge \theta^{-1}_k .
\]
We obtain $G^{i,-2}_{jk} = \frac{\partial A^i_j}{\partial q_k}$ from condition $\overline\Omega^i_j[\overline\omega_x\wedge
\overline\omega_{-1}^k]=0.$

Proceed now to the degree 3.
The f\/irst part of degree 3 we need to compute is $\Omega^i_{-1}$:
\begin{gather*}
 \Omega^i_{-1} \mod \langle \theta_{-2} \wedge \theta_{-3},\theta_{-3} \wedge \theta_{-3}
\rangle   = \frac{\partial f^i}{\partial y_i} \theta_x \wedge \theta^j_{-3} + \frac{\partial B^i_j}{\partial
p_k}\theta^k_{-2} \wedge
\theta^j_{-2}\\
\qquad{}
+ \frac{\partial C^i_j}{\partial x}\theta_x \wedge \theta^j_{-3} + \frac{\partial C^i_j}{\partial q_k}
\theta^k_{-1} \wedge \theta^j_{-3}
-2F^{-3}_j \theta^j_{-3}\wedge \theta^i_{-1}  - 2F^{-2}_j\theta^j_{-2}\wedge
B^i_k\theta^k_{-2}
\\
\qquad {}+G^{i,-3}_{jk}\theta^k_{-3}+G^{i,-2}_{jk}\theta^k_{-2} \wedge B^i_j
\theta^j_{-2}+G^{i,x}_j \theta_x \wedge C^j_k \theta^k_{-3}+2H^{-2}_j\theta^j_{-2}\wedge
\theta^i_{-2} .
\end{gather*}
Wilczynski invariant $W_3$ appears as the $\Omega^i_{-1} [\theta _x \wedge \theta^j_{-3}]$ coef\/f\/icient:
\begin{gather*} \label{11}
\frac{\partial f^i}{\partial y_j} +\frac{d C^i_j}{d x}+A^i_k
C^k_j +2H^x A^i_j.
\end{gather*}
Direct computation shows that:
\begin{gather*}
\Omega^i_{-1} [\theta _x \wedge \theta^j_{-3}]=
\frac{\partial f^i}{\partial y^j} +\frac{1}{3} \frac{\partial f^i}{\partial q^k}\frac{\partial f^k}{\partial p^j} - \frac{d}{dx} \frac{\partial f^i}{\partial
p^j} + \frac{2}{3}\frac{d^2}{dx^2} \frac{\partial f^i}{\partial
q^j} - \frac{2}{27}\left(\frac{\partial f^i}{\partial q^j}\right)^3 \\
\phantom{\Omega^i_{-1} [\theta _x \wedge \theta^j_{-3}]=}{} -\frac{4}{9}
\frac{\partial f^i}{\partial q^k}\frac{d}{dx}\frac{\partial f^k}{\partial q^j}
 -
\frac{2}{9} \frac{d}{dx}\left(\frac{\partial f^i}{\partial q^k}\right)\frac{\partial f^k}{\partial q^j} - 2\delta^i_jH^x .
\end{gather*}
Denote invariant $\Omega^i_{-1} [\theta _x \wedge \theta^j_{-3}]+\frac{1}{2}\frac{d}{dx}W_2$ as $W_3$. Invariant $W_3$ is equivalent to the fundamental invariant $\Omega^i_{-1} [\theta _x \wedge \theta^j_{-3}]$. It means that after replacing $\Omega^i_{-1} [\theta _x \wedge \theta^j_{-3}]$ with $W_3$ the system would remain fundamental. Explicitly the Wilczynski invariant~$W_3$~is:
\begin{gather*}
W_3 = \frac{\partial f^i}{\partial y^j} +\frac{1}{3} \frac{\partial f^i}{\partial q^k}\frac{\partial f^k}{\partial p^j} -\frac{1}{2} \frac{d}{dx} \frac{\partial f^i}{\partial
p^j} + \frac{1}{6}\frac{d^2}{dx^2} \frac{\partial f^i}{\partial
q^j} \\
\phantom{W_3 =}{} - \frac{2}{27}\left(\frac{\partial f^i}{\partial q^j}\right)^3
- \frac{1}{18}
\frac{\partial f^i}{\partial q^k}\frac{d}{dx}\frac{\partial f^k}{\partial q^j}
 -\frac{5}{18} \frac{d}{dx}\left(\frac{\partial f^i}{\partial q^k}\right)\frac{\partial f^k}{\partial q^j} .
\end{gather*}
An  expression $(\frac{\partial f^i}{\partial q^j})^3 $ here is the third power of the matrix $\frac{\partial f^i}{\partial q^j}.$
Note that invariant $W_3$ has known analogue in the case of one dif\/ferential equation of third order:
\[ \frac{\partial f}{\partial y} +\frac{1}{3} \frac{\partial f}{\partial q}\frac{\partial f}{\partial p}- \frac{1}{2} \frac{d}{dx} \frac{\partial f}{\partial
p}+\frac{1}{6}\frac{d^2}{dx^2} \frac{\partial f}{\partial
q} - \frac{2}{27}\left(\frac{\partial f}{\partial q}\right)^3- \frac{1}{3}
\frac{\partial f}{\partial q}\frac{d}{dx}\frac{\partial f^k}{\partial q^j}
 .\]
 The reader can f\/ind this invariant for example in  Chern work \cite{chern}; also see
Sato and Yoshikawa~\cite{satoyo}.

Let us compute the third degree normalization conditions.
   \begin{gather*}
 \Omega_h \mod \langle \theta_{-2} \wedge \theta_{-2},\theta_{-3} \rangle   \\
\qquad{}=F_j^{-3}\theta_x \wedge \theta^j_{-2} + \frac{d F^{-2}_j}{d x} \theta_x \wedge \theta^j_{-2} +
\frac{\partial F^{-2}_j}{\partial q_k} \theta^k_{-1} \wedge \theta^j_{-2}- \theta_x \wedge A^{-2}_j
\theta^j_{-2} .
   \end{gather*}
Thus:
 \[ \Omega_h[\theta_x \wedge \theta^j_{-2}]=-H^{-2}_j + F^{-3}_j+\frac{d F^{-2}_j}{d x} .\]
 Normalizing this coef\/f\/icient to 0 we obtain:
\[  F^{-3}_j =H^{-2}_j -\frac{d F^{-2}_j}{d x}.\]
Next,
 \begin{gather*}
  \Omega^i_j \mod \langle \theta_{-2} \wedge \theta_{-2},\theta_{-3} \rangle =   \frac{\partial
A^i_j}{\partial p_k} \theta^k_{-2} \wedge \theta_x + \frac{d G^{i,-2}_{jk}}{dx} \theta_x \wedge
\theta^k_{-2} +G^{i,-3}_{jk}\theta_x \wedge \theta^k_{-2} \\
\qquad{}+G^{i,x}_k \theta_x \wedge G^{k,-2}_{jl}\theta^l_{-2}+ G^{i,-2}_{kl}\theta^{k}_{-2} \wedge G^{l,x}_j \theta_x .
 \end{gather*}
We have:
\[
\Omega^i_j\big[ \theta_x \wedge \theta^k_{-2} \big]=-\frac{\partial A^i_j}{\partial p_k}
+ \frac{d G^{i,-2}_{jk}}{dx} +G^{i,-3}_{jk} +G^{i,x}_l G^{l,-2}_{jk} - G^{i,-2}_{lk}  G^{l,x}_j .
\]
Assuming this coef\/f\/icient is equal to 0 we get:
\[
G^{i,-3}_{jk} = \frac{\partial A^i_j}{\partial p_k}
- \frac{d G^{i,-2}_{jk}}{dx} - G^{i,x}_l G^{l,-2}_{jk} + G^{i,-2}_{lk}  G^{lx}_j .
\]
Finally,
 \begin{gather*}
  \Omega_y \mod \langle \theta_{-2} ,\theta_{-3} \rangle =   \frac{\partial
H^x}{\partial q_i} \theta^j_{-1} \wedge \theta_x + \frac{d H^{-1}_j}{d x} \theta_x \wedge
\theta^j_{-1}\\
\qquad{} + \frac{\partial H^{-1}_j}{\partial q_k} \theta^k_{-1} \wedge
\theta^j_{-1} +H^{-1}_j \frac{\partial f^j}{\partial q_k} \theta_x \wedge
\theta^k_{-1} + H^{-2}_j \theta_x \wedge \theta^j_{-1} .
 \end{gather*}
The coef\/f\/icient  $\Omega_y[\theta^j_{-1} \wedge \theta_x]$ is the following:
\[ \frac{\partial H^x}{\partial q_j} - \frac{d H^{-1}_j}{d x} -
H^{-1}_k  \frac{\partial f^k}{\partial q_j} - H^{-2}_j .\]
Normalizing it to 0 we obtain:
\[ H^{-2}_j = \frac{\partial H^x}{\partial q_j} - \frac{d H^{-1}_j}{d x} -
H^{-1}_k  \frac{\partial f^k}{\partial q_k} . \]
The last coef\/f\/icient we need in degree 3 is $\Omega_y [ \theta^k_{-1} \wedge \theta^j_{-1}]$:
\[  \frac{\partial
H^{-1}_j}{\partial q_k} - \frac{\partial H^{-1}_k}{\partial q_j} = 0  .\]

In the degree 4 we need to compute only one coef\/f\/icient of curvature:
 \begin{gather*}
\Omega_y \mod \langle \theta_{-2} \wedge \theta_{-2} ,\theta_{-3} \rangle =  \frac{\partial
H^x}{\partial p_j} \theta^j_{-2} \wedge \theta_x + \frac{\partial H^{-1}_j}{\partial p_k} \theta^k_{-2} \wedge
\theta^j_{-1} + \frac{d H^{-2}_j}{d x} \theta_x \wedge
\theta^j_{-2}  \\
{}+ \frac{\partial H^{-2}_j}{\partial q_k} \theta^{-1}_k \wedge
\theta^j_{-2}  + H^{-1}_j \frac{\partial f^j}{\partial p_k} \theta_x \wedge
\theta^k_{-2} + H^{-3}_j \theta_x \wedge \theta^j_{-2} - 2F^{-2}_j \theta^j_{-2} \wedge \big(H^x
\theta_x + H^{-1}_k \theta^k_{-1}\big) .
\end{gather*}
The Cartan connection coef\/f\/icient $\Omega_y[\theta_x \wedge \theta^j_{-2}]$ has the following form:
\[ - \frac{\partial H^x}{\partial p_j} + \frac{d H^{-2}_j}{d x} +
H^{-1}_k  \frac{\partial f^k}{\partial q_j} - H^{-3}_j .\]
Assuming it is equal to 0 we get:
\[ H^{-3}_j = - \frac{\partial H^x}{\partial p_j} + \frac{d H^{-2}_j}{d x} +
H^{-1}_k  \frac{\partial f^k}{\partial q_j} .\]
Finally, invariant $I_4$ is the tensor $\Omega_y[\theta^k_{-1} \wedge \theta^j_{-2}]$:
\[  - \frac{\partial H^{-1}_k}{\partial p_j} + \frac{\partial
H^{-2}_j}{\partial q^k} +
2 H^{-1}_j  H^{-1}_k .\]

\pdfbookmark[1]{References}{ref}
\LastPageEnding


\begin{thebibliography}{99}

\footnotesize\itemsep=0pt

\bibitem{chern}
 Chern S.-S.,
The geometry of the dif\/ferential equation $y''' = F (x, y, y' , y''  )$,
{\em Sci. Rep. Nat. Tsing Hua Univ.~(A)} {\bf 4} (1940), 97--111.

\bibitem{dou}
Doubrov B.,
Contact trivialization of ordinary dif\/ferential equations,
in Dif\/ferential Geometry and Its Applications (Opava, 2001), {\it Math. Publ.}, Vol.~3, Silesian Univ. Opava, Opava, 2001, 73--84.

\bibitem{dkm}
Doubrov B., Komrakov B., Morimoto T.,
Equivalence of holonomic dif\/ferential equations,
{\em Lobachevskii J. Math.} {\bf 3} (1999), 39--71.

\bibitem{fels}
Fels M.,
The equivalence problem for systems of second-order ordinary dif\/ferential equations,
\href{http://dx.doi.org/10.1112/plms/s3-71.1.221}{{\em Proc. London Math. Soc.}}  {\bf 71} (1995), 221--240.

\bibitem{lie}
Lie S.,
Vorlesungen \"uber Dif\/ferentialgleichungen mit bekannten  inf\/initesimalen Transformationen,
Leipzig, Teubner, 1891.

\bibitem{medv}
Medvedev A.,
Geometry of third order ODE systems,
{\em Arch. Math. (Brno)} {\bf 46} (2010), 351--361.

\bibitem{mor93}
Morimoto T.,
Geometric structures on f\/iltered manifolds,
{\em Hokkaido Math.~J.} {\bf 22} (1993), 263--347.

  \bibitem{satoyo}
 Sato H., Yoshikawa A.Y.,
Third order ordinary dif\/ferential equations and Legendre connections,
\href{http://dx.doi.org/10.2969/jmsj/05040993}{{\em J. Math. Soc. Japan}}  {\bf 50} (1998), 993--1013.

\bibitem{tan70}
Tanaka N.,
On dif\/ferential systems, graded Lie algebras and pseudo-groups,
{\em J.~Math.  Kyoto. Univ.} {\bf 10} (1970), 1--82.

\bibitem{tan79}
Tanaka N.,
On the equivalence problems associated with simple graded Lie   algebras,
 {\em Hokkaido Math.~J.} {\bf 8} (1979), 23--84.

\bibitem{tresse}
 Tresse M.A.,
D\'etermination des invariants ponctuels de l'\'equation  dif\/f\'erentielle ordinaire du second ordre $y'' =\omega(x,y,y')$,
 Leipzig, Hirzel, 1896.

\bibitem{yano}
Yano K.,
 The theory of Lie derivatives and its applications,
  North Holland Publishing Co., Amsterdam, 1957.



\end{thebibliography}
\end{document}